\documentstyle{article}
\begin{document}

\begin{center}
\Large{ BRST OPERATOR FOR QUANTUM \\
LIE ALGEBRAS AND DIFFERENTIAL CALCULUS \\
ON QUANTUM GROUPS}
\end{center}

\vspace{.5cm}

\begin{center}
\large{ A.P.\,ISAEV${}^{a \,1}$ and O.\,OGIEVETSKY${}^{b \, 2}$ }
\end{center}

\begin{center}
${}^{a}$ Bogoliubov Theoretical Laboratory,
Joint Institute for Nuclear Research,
Dubna, Moscow region 141980, Russia. \\
${}^{b}$ Center of Theoretical Physics, Luminy,
13288 Marseille, France \\ 
and P. N. Lebedev Physical Institute, Theoretical Department,
Leninsky pr. 53, 117924 Moscow, Russia.
\end{center}

\vspace{1cm}

\begin{center}
ABSTRACT
\end{center}

For a Hopf algebra ${\cal A}$,
we define the structures of 
differential complexes on two dual exterior Hopf algebras: 1) an exterior
extension of ${\cal A}$  and 
2) an exterior extension of the dual algebra ${\cal A}^*$. 
The Heisenberg double of these two exterior Hopf algebras
defines the differential algebra for
the Cartan differential calculus on ${\cal A}$.
The first differential complex is an analog of the de Rham complex.
In the situation when ${\cal A}^*$ is a universal
enveloping of a Lie (super)algebra
the second complex coincides with the standard complex.
The differential is realized as an (anti)commutator with
a BRST- operator $Q$. A recurrent relation which defines uniquely
the operator $Q$ is given.
The BRST and anti-BRST operators  are constructed explicitly
and the Hodge decomposition theorem
is formulated for the case of the quantum Lie algebra $U_q(gl(N))$.

\vspace{1cm}

${}^1$ E-mail: isaevap@thsun1.jinr.ru

\vspace{0.5cm}
${}^2$ E-mail: oleg@cpt.univ-mrs.fr

\pagebreak

\section{Introduction}

The theory of the bicovariant differential calculus on a Hopf algebra
${\cal A}$ has been developed by S.Woronowicz \cite{Wor} 
(for reviews see \cite{Cast}, \cite{Schm}, \cite{Is}) on the
basis of the axiomatics for the noncommutative differential geometry
proposed by A.Connes \cite{Con}. 
Then it was shown in \cite{Man}, \cite{Jur} that the Woronowicz
theory could be applied to the description of differential calculi on
quantum analogs of some Lie groups and can be adopted to the
$R$-matrix formalism of \cite{FRT}. 
The special example of $Fun(GL_q(N))$ has been considered
in detail in \cite{Man}, \cite{SIP}, \cite{SWZ}, \cite{IsPya} 
(see also \cite{Is} and references therein). 

In \cite{Man}, \cite{IsPya}, \cite{Is} it was established that 
the $GL_q(N)$- differential algebra is a Poincar\'{e}-Birkhoff-Witt (PBW) 
type algebra
and therefore indeed describes the quantum deformation of the classical
differential algebras over $GL(N)$. 
The $SL_q(N)$- differential algebra of PBW type has been
constructed in \cite{FP} (see also  \cite{AAM}). 

However, further
analysis \cite{AIP} has shown that the bicovariant differential
algebras \cite{Jur}, \cite{{CWSWW}} of the Woronowicz calculus on the 
$Fun(SO_q(N))$ and $Fun(Sp_q(2n))$ (actually their bicovariant subalgebras of
differential 1-forms) have incorrect dimensions and are not of PBW
type. Thus, these noncommutative algebras could not be interpreted 
as the quantum deformations of the corresponding classical differential
algebras. 

In this paper we relate the Woronowicz theory 
to a deformation of the BRST theory
(for review of the BRST theory see \cite{HT}; the applications
of the BRST theory to the Lie algebra
cohomology theory can be found in \cite{vanH}, \cite{CAMP}). This
relation has already been exploited in \cite{Ber}, \cite{Wat} 
in the context
of the discussion of the so-called
"quantum group gauge theories" (about the $q$-group gauge theories 
and related noncommutative geometries see
e.g. \cite{AVIP} and references therein).  The main idea is
that the Woronowicz exterior differential map $d$ acting on the 
exterior extension of the Hopf
algebra ${\cal A}$ should be generated by
a nilpotent operator which is nothing but the BRST charge
related to the deformed algebra of the vector fields over ${\cal A}$. 
The construction
of the explicit formula for this BRST operator is the main result of
the present paper. 

In the first Section we explain the notion of the quantum Lie algebra.
In the second Section we collect all results about the Cartan differential
calculus on Hopf algebras (this is an extension of the Woronowicz
differential calculus) and show how the quantum Lie algebras appear
naturally in the context of these calculi. The BRST operator
for the quantum Lie algebras is constructed in Section 3. 
In Section 4 we consider the special case
of the quantum Lie algebra $U_q(gl(N))$ in detail.
The BRST and anti-BRST operators  are constructed explicitly
and the Hodge decomposition theorem
is formulated for this case.

\section{Quantum Lie Algebras}

A quantum Lie algebra \cite{Wor1}, \cite{Wor},
\cite{Ber}, \cite{Cast} is defined by two tensors
$C^k_{ij}$ and $\sigma^{mk}_{ij}$ (indices belong
to some set ${\cal N}$, say, ${\cal N}= \{ 1, \dots, N \}$). 
By definition, the matrix $\sigma^{mk}_{ij}$ has an eigenvalue
1; one demands that $(P_{(1)})^{mk}_{ij} \, C^n_{mk} = 0$, where
$P_{(1)}$ is a projector on the eigenspace of $\sigma$ corresponding to the
eigenvalue 1. 

By definition, a quantum Lie algebra $\Gamma$
is generated by elements $\chi_i$, $i=1, \dots, N$, subjected
to relations
 \begin{equation}
 \label{int1}
\chi_i \, \chi_j - \sigma^{mk}_{ij} \chi_m \, \chi_k = C^k_{ij} \,
\chi_k \; .
 \end{equation}
Here the structure constants $C^k_{ij}$ obey
 \begin{equation}
 \label{int2}
 \begin{array}{c}
C^p_{ni} \, C^{l}_{pj} = \sigma^{mk}_{ij} \, C^p_{nm} \,
C^{l}_{pk} + C^p_{ij} \, C^{l}_{np} \; \Leftrightarrow \\ \\
\Leftrightarrow \;
C^{<1|}_{|1 2>} \, C^{<4|}_{|1 3>} = \sigma_{23} \, C^{<1|}_{|1 2>} \,
C^{<4|}_{|1 3>} + C^{<3|}_{|23>} \, C^{<4|}_{|1 3>} \; ,
\end{array}
 \end{equation}
 \begin{equation}
 \label{int3}
C^k_{ni} \, \sigma^{pm}_{kq} =  \sigma^{sj}_{iq} \,  
\sigma^{pk}_{ns} \, C^m_{kj} 
\;\;\; \Leftrightarrow \;\;\;
C^{<1|}_{|1 2>} \, \sigma_{13} = 
\sigma_{23} \, \sigma_{12} \, C^{<3|}_{|2 3>} \; ,
 \end{equation}
 
 \begin{equation}
 \label{int3a}
 \begin{array}{c}
(\sigma^{pj}_{im} \, C^n_{qp} + \delta^n_q \, C^j_{im} ) \,
\sigma^{ks}_{nj}
= \sigma^{jn}_{qi} \, (\sigma^{ps}_{nm} \, C^k_{jp} +
\delta^k_j  \, C^{s}_{nm})  \;\;\; \Leftrightarrow  \\ \\
 \Leftrightarrow \;\;\;
(\sigma_{23} \, C^{<1|}_{|12>} +  C^{<3|}_{|23>} ) \,
\sigma_{13} = \sigma_{12} \, (\sigma_{23} \, C^{<1|}_{|12>} +
C^{<3|}_{|23>} ) \; .
\end{array}
 \end{equation}
The matrix $\sigma^{mk}_{ij}$ satisfies the Yang-Baxter equation
 \begin{equation}
 \label{int1a}
\sigma^{j_1 j_2}_{i_1 i_2} \, \sigma^{n_2 k_3}_{j_2 i_3} \,
\sigma^{k_1 k_2}_{j_1 n_2} = \sigma^{j_2 j_3}_{i_2 i_3} \,
\sigma^{k_1 n_2}_{i_1 j_2} \, \sigma^{k_2 k_3}_{n_2 j_3} \;\;\;
\Leftrightarrow \;\;\;  \sigma_{12} \, \sigma_{23} \, \sigma_{12} =
\sigma_{23} \, \sigma_{12} \, \sigma_{23} \; .
 \end{equation}
In the right hand 
side of (\ref{int3})-(\ref{int1a}) we
use the FRT matrix notations \cite{FRT}; indices $\{ 1,2,3, \dots \}$
are the numbers of vector spaces, 
{\it e.g.}, $f_1 := f^{i_1}_{j_1}$
is a matrix which acts in the first vector space.
Additionally, we use  
incoming and outcoming indices, {\it e.g.},
$\Omega^{<1|} := \Omega^{i_1}$ and  $\gamma_{|1>} := \gamma_{j_1}$
denote a covector with one outcoming index and 
a vector with one incoming index respectively. 
Thus, in this notation, the matrix $f_1$ can be written as
$f_1 = f^{<1|}_{|1>}$.

Note that relations (\ref{int2}) - (\ref{int1a}) can be justified
if we consider a monomial 
$\chi_{|1>} \chi_{|2>} \chi_{|3>}$ 
of degree three
and reorder it in two different
ways using the defining relations (\ref{int1}). 
A demand that the results coincide implies 
\begin{equation}
\label{co1}
\begin{array}{c}
0= (\sigma_{12} \, \sigma_{23} \, \sigma_{12}-
\sigma_{23} \, \sigma_{12} \, \sigma_{23})
\, \chi_{|1>} \, \chi_{|2>} \, \chi_{|3>} \\
+ \left[ \left( \sigma_{23} \,  C^{<1|}_{|12>}  +
C^{<3|}_{|23>} \right) \, \sigma_{13}
- \sigma_{12} \, \left( \sigma_{23} \, C^{<1|}_{|12>} 
+C^{<3|}_{|23>} \right) \right] \, \chi_{|1>} \chi_{|3>} \\
-\left[ C^{<1|}_{|12>} \, \sigma_{13} -
\sigma_{23}  \, \sigma_{12} \, C^{<3|}_{|23>} \right] \, \chi_{|1>} 
\, \chi_{|3>} \\
- \left( (1  - \sigma_{23}) \, C^{<1|}_{|12>}  +
C^{<3|}_{|23>} \right) \, C^{<4|}_{|13>} \, \chi_{|4>}.
\end{array}
\end{equation}
This is indeed an identity: the cubic term vanishes because of 
(\ref{int1a}), the quadratic terms vanish
in view of (\ref{int3a}) and (\ref{int3}) while the last term 
disappears due to the Jacobi identity (\ref{int2}). 

\vspace{.1cm}
\noindent
{\bf Remark 1.} Note that to obtain the identity (\ref{co1})
it is not necessary to require precisely the relations
(\ref{int3}) and (\ref{int3a}), only their combination enters
in (\ref{co1}). To obtain relations
(\ref{int3}) and (\ref{int3a}) as consistency conditions, 
one has to consider the following situation. For each integer
$M$, let $\chi^{(M)}_i$ be a copy of the generators $\chi_i$.
Assume that the relations between different copies are given by
\begin{equation} 
\chi^{(K)}_{|1>}\chi^{(M)}_{|2>} =
\sigma_{12}\chi^{(M)}_{|1>}\chi^{(K)}_{|2>} +C_{|12>}^{<3|}
\chi^{(K)}_{|3>}\ \ \ \ {\mathrm for}\ \ \ \ K<M\ .
\end{equation}
Then, ordering in two different ways an expression
$\chi^{(L)}_{|1>}\chi^{(M)}_{|2>}\chi^{(K)}_{|3>}$ with
$L<M<K$, one obtains (\ref{int1a}) in terms cubic in $\chi$, and,
in lower odred terms, an identity
\begin{equation}
\begin{array}{l}
\left[ \left( \sigma_{23} \,  C^{<1|}_{|12>}  +
C^{<3|}_{|23>} \right) \, \sigma_{13}
- \sigma_{12} \, \left( \sigma_{23} \, C^{<1|}_{|12>}  +
C^{<3|}_{|23>} \right) \right] \, \chi^{(M)}_{|1>} \chi^{(L)}_{|3>} \\
-\left[ C^{<1|}_{|12>} \, \sigma_{13} -
\sigma_{23} \, \sigma_{12} \, C^{<3|}_{|23>} \right] \, 
\chi^{(K)}_{|1>} \chi^{(L)}_{|3>} \\
- \left( (1  - \sigma_{23}) \,  C^{<1|}_{|12>}  -
C^{<3|}_{|23>} \right) \, C^{<4|}_{|13>} \, \chi^{(L)}_{|4>} =0\, .
\end{array}\end{equation}
This identity is completely equivalent to the set of the
relations 
(\ref{int2}), (\ref{int3}) and (\ref{int3a}).

\vspace{.1cm}
\noindent
{\bf Remark 2.} The Jacobi identity (\ref{int2}) implies the existence of
the adjoint representation, in which the generator $\chi_i$ is represented
by a matrix $({\mathrm ad}({\chi_i}))^j_k=C^j_{ki}$.
Also, any quantum Lie algebra possesses a trivial one dimensional
representation, in which each generator $\chi_i$ acts as zero.

\vspace{.1cm}
\noindent
{\bf Remark 3.}
The quantum Lie algebras defined by equations 
(\ref{int1})-(\ref{int1a}) generalize the
usual Lie (super-)algebras. In the non-deformed case, when
$\sigma^{mk}_{ij} = (-1)^{(m)(k)} \, \delta^m_j \, \delta^k_i$
is a super-permutation matrix (here $\sigma^2 =1$
and (\ref{int1a}) is fulfilled;
$(m)=0,1$ $mod(2)$ is the parity of a generator $\chi_m$), equations
(\ref{int1}) and (\ref{int2}) coincide with the defining 
relations and the Jacobi identities for Lie (super)-algebras.
Equation (\ref{int3}) is then equivalent to the 
$Z_2$-homogeneity condition $C^i_{jk} =0$ for $(i) \neq (j)+(k)$.
Equation (\ref{int3a}) follows from (\ref{int3}).

\vspace{.1cm}
\noindent
{\bf Remark 4.} It was noted in \cite{Ber} that the relations 
(\ref{int2}) - (\ref{int1a})  can be encoded as 
the Yang-Baxter equation for the matrix
$S^{AB}_{CD}$ where capital latin indices $A$, $B$, ... belong to
the set $0\cup {\cal N}$ (for the matrix $\sigma^{ij}_{kl}$
the indices belong to the set ${\cal N}$). The matrix $S$ is defined
by 
\begin{equation} 
S^{ij}_{kl}=\sigma^{ij}_{kl}\ ,\ 
S^{0j}_{kl}=C^j_{kl}\ ,\ S^{0A}_{0B}=\delta^A_B 
\end{equation}
and the other components of $S$ are zeros.

\section{Cartan Differential Calculus on Hopf Algebras}
In this Section we explain that the quantum Lie algebras
(\ref{int1}) appear naturally (as quantum analogs of vector fields)
in the context of the bicovariant differential calculus on the
Hopf algebras.

\subsection{Exterior Hopf algebras}
Let ${\cal A}$ $(\Delta, \epsilon , S)$ be a Hopf algebra and
${\cal A}^*$ is a Hopf dual to ${\cal A}$. 
The comultiplication and left-right ${\cal A}$- coactions on
${\cal A}$ are 
 \begin{equation}
 \label{1}
\Delta (a) = a_{(1)} \otimes a_{(2)} \equiv \Delta_{L,R} (a) \; ,
\;\;\; a \in {\cal A} \; .
 \end{equation}
Here the Sweedler notation is used. 
One can define \cite{Wor} the bicovariant bimodule $\Gamma^{(1)}$ over
${\cal A}$ as a linear space with left-invariant 
basic elements $\{ \omega^i \}$
such that left and right ${\cal A}$- coactions on $\Gamma^{(1)}$
have the form
 \begin{equation}
 \label{2}
\Delta_L (\omega^i ) = 1 \otimes \omega^i \; , \;\;\;
\Delta_R (\omega^i ) =  \omega^j  \otimes r^i_j
 \end{equation}
where $r^i_j \in {\cal A}$. Since any left-invariant element can be written
as a linear combination of 
the basis elements $\{ \omega^i \}$, one has
$S(a_{(1)}) \, \omega^i \, a_{(2)} = f^{i}_j (a) \, \omega^j$
which is equivalent to the commutation relations of elements
$a \in {\cal A}$ with $\omega^i \in \Gamma^{(1)}$:
 \begin{equation}
 \label{3}
\omega^i \, a = ( f^i_j \triangleright a ) \omega^j = a_{(1)}
f^{i}_j (a_{(2)}) \, \omega^j \; ,
 \end{equation}
where  $f^i_j \in {\cal A}^*$. Covariance of (\ref{3}) under right
${\cal A}$- coaction (\ref{1}), (\ref{2}) requires the main condition
on the elements $f, \, r$:
 \begin{equation}
 \label{4}
(f^j_i \triangleright a ) \, r^i_k = r^j_i \, (a \triangleleft
f^i_k) \; , \;\;\; \forall a \in {\cal A} \; .
 \end{equation}

As it was shown in \cite{Brz} (for more details see \cite{Schm}),
one can construct the exterior Hopf algebra $\Gamma^{\wedge}
= \{ {\cal A} \oplus \Gamma^{(1)} \oplus \Gamma^{(2)} \oplus \dots
\}$ via the Woronowicz's definition of the covariant exterior
product 
 \begin{equation}
 \label{5}
 \omega^i \wedge \omega^j = \omega^i \otimes \omega^j -
\omega^k \otimes \omega^l \, \sigma^{ij}_{kl} \Rightarrow
\omega^{<1|} \wedge \omega^{<2|} = \omega^{<1|} \otimes
\omega^{<2|} \, (1 - \sigma)_{12} \; .
 \end{equation}
where the matrix $\sigma^{ij}_{kl} = f^i_l(r^j_k)$ 
(an analogue of the permutation matrix) satisfy the Yang-Baxter
equation (\ref{int1a}) which follows from (\ref{4}) if we put $a = r^m_n$ and
then take the pairing with $f^l_p$ (here we need to use the explicit
forms for the comultiplications $\Delta(f^i_j)$ and $\Delta(r^i_j)$,
see below).

A generalization of the definition of the wedge product (\ref{5})
to the case of the product of $n$ 1-forms is straightforward
 \begin{equation}
 \label{owedge}
\omega^{<1|} \wedge \omega^{<2|} \wedge \dots \wedge \omega^{<n|}
= \omega^{<1|} \otimes \omega^{<2|} \otimes \dots \otimes
\omega^{<n|} \, A_{1 \rightarrow n} \; .
 \end{equation}
Here the matrix operator $A_{1 \rightarrow n}$ - is an analog of
the antisymmetrizer of $n$ spaces. This operator is defined
inductively (see e.g. \cite{Gur})
 \begin{equation}
 \label{anti}
\! A_{1\rightarrow n} = \left(1 - \sum_{k=1}^{n-1} \, (-1)^{n-k-1} \,
\sigma_{k \rightarrow n} \right)  A_{1\rightarrow n-1} = \left(1
- \sum_{k=2}^{n} \, (-1)^{k} \, \sigma_{k \leftarrow 1}
\right) A_{2\rightarrow n} \; ,
 \end{equation}
where $\sigma_{n \leftarrow k} = \sigma_{n-1 n} \dots \sigma_{k+1 k+2} \, 
\sigma_{k k+1}$, $\sigma_{k \rightarrow n} = \sigma_{k k+1} \, 
\sigma_{k+1 k+2} \dots \sigma_{n-1 n}$ ($n >k$). 
If the sequence of operators $A_{1 \rightarrow n}$ vanishes
at the step $n = h+1$ ($A_{1 \rightarrow n} = 0$ $\forall n >h$)
then the number $h$ is called the height of 
the matrix $\sigma_{12}$.

The extension of the comultiplication (\ref{1}) to 
the whole exterior algebra $\Gamma^{\wedge}$ is defined by
\cite{Brz}
 \begin{equation}
 \label{6}
\Delta (\omega^i) \equiv \Delta_L (\omega^i) + \Delta_R (\omega^i) =
 1 \otimes \omega^i +  \omega^j  \otimes r^i_j \; ,
 \end{equation}
\begin{equation}
 \label{7}
\Delta (\omega \, a) = \Delta (\omega ) \, \Delta (a) \; , \;\;\;
\Delta (a \, \omega) = \Delta (a) \, \Delta (\omega) \; ,
 \end{equation}
where $\otimes$ is a graded tensor product
and the grading is: $deg(\omega)=n$ for $\omega \in \Gamma^{(n)}$.
Associativity condition for $(\omega \, a \, b)$ 
with respect to (\ref{3}) and coassociativity condition for (\ref{6}) 
yield the form of
comultiplications for $f^i_j$ and $r^i_j$
 \begin{equation}
  \label{7a}
 \Delta(f^i_j) = f^i_k \otimes f^k_j \; , \;\;\;
 \Delta(r^i_j) = r^k_j \otimes r^i_k  \; .
 \end{equation}
The other structure mappings for $f$, $r$ and $\omega$'s are
obtained from (\ref{6}), (\ref{7a})
 $$
\epsilon(r^i_j) = \delta^i_j = \epsilon(f^i_j) \; , \;\;\;
S(f^i_k) \,f^k_j = \delta^i_j = S(r^k_j) \, r^i_k \; , 
$$
\begin{equation}
\label{eso}
\epsilon(\omega^i) = 0 \; , \;\;\; S(\omega^i) = - \omega^j \,
S(r^i_j) \; .
\end{equation}

To summarize this subsection we stress that the knowledge of two sets
of elements $\{ r^i_j \} \in {\cal A}$ and $\{f^i_j \} \in {\cal
A}^*$ which satisfy (\ref{4}) is enough to construct a bicovariant
bimodule over ${\cal A}$ and then extend ${\cal A}$ to the
exterior Hopf algebra $\Gamma^{\wedge}$.

\vspace{0.1cm}
\noindent
{\bf Remark.}
The space $\Gamma^{(k)}$ is a subspace in $\omega^{\otimes \, k}$
spanned by tensors $a$ of the form 
\begin{equation}
\label{def}
a = \omega^{i_1} \otimes \dots \otimes \omega^{i_k} \, 
A^{j_1 \dots j_k}_{i_1 \dots i_k} \, a_{j_1 \dots j_k}
\end{equation}
The formula $A_{1 \dots k} f_1 \dots f_k = f_1 \dots f_k \, A_{1 \dots k}$
implies that multiplication of the
elements of $\Gamma^k$ and ${\cal A}$ is compatible
with (\ref{3}).
Given two forms $a \in \Gamma^{(k)}$ (as in (\ref{def})) and 
$b \in \Gamma^{(l)}$ (with coefficients $b_{j_1 \dots j_l}$) 
define their wedge product $a \wedge b \in \Gamma^{(k+l)}$ to be
\begin{equation}
\label{ab}
a \wedge b := \omega^{i_1} \otimes \dots \otimes \omega^{i_{k+l}} \, 
A^{j_1 \dots j_{k+l}}_{i_1 \dots i_{k+l}} \, c_{j_1 \dots j_{k+l}}
\end{equation}
where
\begin{equation}
\label{c}
c_{i_1 \dots i_k j_1 \dots j_l} = 
a_{i_1 \dots i_k} \, b_{j_1 \dots j_l} \; .
\end{equation}
The definition (\ref{ab}), (\ref{c})
of the wedge product of the differential forms is selfconsistent. 
Indeed, one can change $a$ by adding $\delta a$, such that
$A_{1 \dots k} \, \delta a = 0$. Then, 
$A_{1 \dots k+l} \, \delta a \wedge b = 0$, since $A_{1 \dots k+l}$ (in view
of (\ref{anti})) is proportional to $A_{1 \dots k}$,
$A_{1 \dots k+l} = Z \cdot A_{1 \dots k}$ for some $Z$.
It is straightforward to see that this tensor product is associative.

\subsection{Dual exterior Hopf algebra}

By analogy with the construction of the previous subsection one can
define the bicovariant bimodule $\Gamma^{(1)*}$ and exterior Hopf
algebra over the dual algebra ${\cal A}^*$ ($\Delta (h) = h_{(1)}
\otimes h_{(2)}$, $h \in {\cal A}^*$) by introducing two sets of
elements $\{ \overline{r}^i_j \} \in {\cal A}$ and
$\{\overline{f}^i_j \} \in {\cal A}^*$ such that
(cf. with (\ref{3}), (\ref{6}), (\ref{eso}))
\begin{equation}
\label{9} 
\gamma_i \, h = (\overline{r}^j_i \triangleright h) \, \gamma_j 
= h_{(1)} \, < h_{(2)} \, , \, \overline{r}^j_i > \, \gamma_j \; ,
\end{equation}
\begin{equation}
\label{8}
\Delta(\gamma_i) = 1 \otimes \gamma_i + \gamma_j \otimes
\overline{f}^j_i \; ,
\end{equation}
\begin{equation}
\label{8g}
\epsilon (\gamma_i) = 0 \; , \;\;\; S(\gamma_i) = - \gamma_j \,
S(\overline{f}^j_i) \; \Rightarrow \; S^{-1}(\gamma_i) = -
S^{-1}(\overline{f}^j_i) \, \gamma_j \; ,
\end{equation}
where $<h, \, a> = h(a)$ is a dual pairing for Hopf algebras
${\cal A}$, ${\cal A}^*$ and
the elements $\{ \gamma_i \}$ form the left-invariant basis
for the bimodule $\Gamma^{(1)*}$. As above, we have from $\Delta
(\gamma_i \, h) = \Delta (\gamma_i) \, \Delta (h) $ the relation
(cf. with (\ref{4}))
\begin{equation}
\label{10}
(\overline{r}^j_i \triangleright h ) \, \overline{f}^k_j =
\overline{f}^j_i \, (h \triangleleft \overline{r}^k_j ) \; , \;\;\; 
\forall h \in {\cal A}^* \; .
\end{equation}
and the wedge product in the exterior Hopf algebra $\Gamma^{\wedge *}
= {\cal A}^* \oplus \Gamma^{(1)*} \oplus \Gamma^{(2)*} \oplus
\dots$ ($deg(\gamma)=-n$ for $\gamma \in \Gamma^{(n)*}$)
is defined by
\begin{equation}
\label{11}
\gamma_i \wedge \gamma_j = \gamma_i \otimes \gamma_j -
\overline{\sigma}^{kl}_{ij} \,  \gamma_k \otimes \gamma_l \; ,
\;\;\; \overline{\sigma}^{kl}_{ij} := \overline{r}^l_i
(\overline{f}^k_j) \; .
\end{equation}

Let us consider the special case when $dim (\Gamma^{(1)}) = dim
(\Gamma^{(1)*})$ and formulate the conditions when the exterior
Hopf algebras $\Gamma^{\wedge}$ and $\Gamma^{\wedge *}$ are Hopf
dual to each other. One can extend the pairing for the algebras
${\cal A}$ and ${\cal A}^*$ to
the non-degenerate pairing of the
algebras $\Gamma^{\wedge}$ and $\Gamma^{\wedge *}$:
\begin{equation}
\label{12}
< \gamma_i , \, \omega^j > = \delta^j_i \; , \; \; \;
<\Gamma^{(n) *} , \, \Gamma^{(m)} > \sim \delta^{nm} \; .
\end{equation}
This pairing is compatible with the grading. The 
relations (\ref{6}), (\ref{8}) and (\ref{12}) give
\begin{equation}
\label{13}
< h \, \gamma_i , \, a \, \omega^j > =
<h, \, a> \, \delta^j_i \; , \;\;\; h \in {\cal A}^* , \;\;
a \in {\cal A} \; .
\end{equation}
Now it is clear that the equations (\ref{12}) and (\ref{13}) relate the
sets of elements $\{ r, \, f \}$ with $\{ \overline{r}, \,
\overline{f} \}$. Indeed, from one side we have $ < \gamma_i , \,
\omega^k \, a> = < \Delta(\gamma_i) , \, \omega^k \otimes a> =
\overline{f}^k_i(a)$, $\forall a$ but on the other side we deduce $ <
\gamma_i , \, \omega^k \, a> = < \gamma_i , \, (f^k_j
\triangleright a) \, \omega^j > = f^k_i(a) $ and therefore
$\overline{f}^k_i = f^k_i$. Considering the
pairing $< \gamma_i \, h , \omega^k >$, we obtain $\overline{r}^k_i
= r^k_i$ and, thus, $\overline{\sigma}^{ij}_{kl} =
\sigma^{ij}_{kl}$. It leads to the definition of wedge product of
$n$ elements $\gamma_i$ (see (\ref{11})):
\begin{equation}
\label{13a}
\gamma_{|1>} \wedge \gamma_{|2>} \dots \wedge \gamma_{|n>} =
A_{1 \rightarrow n} \, \gamma_{|1>} \otimes \gamma_{|2>} \dots
\otimes \gamma_{|n>} \; ,
\end{equation}
where operators $A_{1 \rightarrow n}$ are the same as in
(\ref{anti}).

\subsection{Differential $d$ and the algebra of vector fields}

Suppose that there exists a differential map $d$: $\Gamma^{(n)}
\rightarrow \Gamma^{(n+1)}$ $(\Gamma^{(0)} = {\cal A})$
which squares to 0 and satisfies the Leibniz rule
$(\omega, \, \omega_i \in \Gamma^{\wedge})$ 
\begin{equation}
\label{14}
d^2 (\omega) = 0 \; , \;\;\;
d(\omega_1 \, \omega_2 ) = d(\omega_1) \, \omega_2 +
(-1)^{deg (\omega_1)} \, \omega_1 \, d(\omega_2) \; .
\end{equation}
As it was shown in \cite{Wor} the left 
and right coactions: $\Delta_L(d(a)) =
a_{(1)} \otimes d(a_{(2)})$, $\Delta_R(d(a)) = d(a_{(1)}) \otimes
a_{(2)}$ are compatible with the definition of the bicovariant
bimodule $\Gamma^{(1)}$. From (\ref{6}) we obtain the "Leibniz
rule" for the comultiplication
\begin{equation}
\label{14b}
\Delta(d(a)) = (d \otimes id + id \otimes d) \Delta(a) \Rightarrow
\Delta \, d = (d \otimes id + id \otimes d) \Delta \; ,
\end{equation}
which can be extended to the whole exterior algebra
$\Gamma^{\wedge}$ (taking into account that $\otimes$
is the graded tensor product). Since the 1-form $\omega = S(a_{(1)})
\, da_{(2)} \in \Gamma^{(1)}$ $(\forall a \in {\cal A})$ 
is left-invariant, it has to be
expanded over the left-invariant basis $\{ \omega^i \}$:
\begin{equation}
\label{15}
S(a_{(1)}) \, da_{(2)} = \sum_i \, \chi_i(a) \, \omega^i \;
\Rightarrow \; d(a) = (\chi_i \triangleright a) \, \omega^i \; ,
\end{equation}
where $\chi_i(a)$ are some coefficients
and $\chi_i \in {\cal A}^*$. Applying the
Leibniz rules (\ref{14}) and (\ref{14b}) to eq. (\ref{15}) gives
\begin{equation}
\label{16dd}
\Delta(\chi_i) =
\chi_j \otimes f^j_i + 1 \otimes \chi_i  \; , 
\end{equation}
\begin{equation}
\label{16}
a \triangleleft \chi_i = (\chi_j \triangleright a) r^j_i \ .
\end{equation} 
Equations (\ref{15}) and (\ref{16dd}) lead to the definition of
the antipode and the counit for the elements $\chi_i$: 
\begin{equation}
\label{16a}
\chi_i(I) = \epsilon(\chi_i) = 0 \; , \;\;\;
S(\chi_i) = - \chi_j \, S(f^j_i) \Rightarrow
\chi_i = - S(\chi_j) \, f^j_i \; .
\end{equation}
According to \cite{Wor}, $\chi_i$ are interpreted  as  vector
fields over the Hopf algebra ${\cal A}$.

One can obtain the commutation rules for
elements $\chi_i$ with arbitrary $h \in {\cal A}^*$:
\begin{equation}
\label{19}
\chi_i \, h = (r_i^k \triangleright h ) \, \chi_k  \; ,
\end{equation}
from the requirement that the pairing of ${\cal A}$
and  ${\cal A}^*$ is non-degenerate.
Indeed, we have from (\ref{16}) 
$$
< \chi_i \, h , \, a > = < h , \,
a \triangleleft \chi_i > = < h , \, (\chi_k \triangleright a)
r^k_i > = < (r^k_i \triangleright h) \, \chi_k , \, a >
$$ 
and, since $a \in {\cal A}$ is an arbitrary element, we deduce
(\ref{19}). Equations (\ref{19}) for $h = \chi_j$ give the
defining relations (\ref{int1}) for the elements $\chi_i$:
\begin{equation}
\label{21}
\chi_i \, \chi_j =
(r_i^k \triangleright \chi_j ) \, \chi_k  =
\sigma^{mk}_{ij} \chi_m \, \chi_k
+ C^k_{ij} \, \chi_k  \; , \;\;\;
C^k_{ij} = <\chi_j , \, r_i^k > \; .
\end{equation} 
The Jacobi
identities (\ref{int2}) for the structure constants $C^i_{jk}$
can be obtained by pairing of eq. (\ref{21}) with $r^l_n$.
Another application of (\ref{16}) is that by taking $a=r^i_j$ one
deduces the relation
\begin{equation}
\label{36}
C^k_{ni} \, r^m_k = r^k_n \, r^j_i \, C^m_{kj} \;
\Leftrightarrow \; C^{<3|}_{|12>} \, r_{3} = r_1 \, r_2 \,
C^{<3|}_{|12>} \;  ,
\end{equation}
which gives the conditions (\ref{int3}) by pairing with $f^p_q$.
Eqs. (\ref{36}) are the invariance condition for
the structure constants $C^k_{ni}$ with respect to the 
rotations by matrices $r^i_j$.

Eqs. (\ref{10}), where $\overline{r}^i_j = r^i_j$ and
$\overline{f}^i_j = f^i_j$, give for $h = f^m_n$ and $h = \chi_n$:
\begin{equation}
\label{sff}
\sigma_{12} \, f_1 \, f_2 = f_1 \, f_2 \, \sigma_{12} \; ,
\end{equation}
\begin{equation}
\label{com2b}
(\sigma^{pj}_{im} \, \chi_p + C^j_{im} ) \, f^k_j = f^j_i \, f^p_m
\, C^k_{jp} + f_i^k \, \chi_m \; ,
\end{equation}
and pairing (\ref{sff}) and (\ref{com2b}) with $r^s_q$ reproduces 
(\ref{int1a}) and (\ref{int3a}).

Now we introduce the set of elements $a^i \in {\cal A}$ such that $<
\chi_i , \, a^j > = \delta^j_i$. By definition we have
$S(a^j_{(1)}) d (a^j_{(2)}) = \omega^j$ and $\epsilon(a^i) = 0$
\cite{Wor}. 

The Maurer-Cartan equation reflects the fact that $d^2(a)=0$ and
this equation can be deduced as follows
\begin{equation}
\label{24}
\begin{array}{c}
d \omega^k = d (S(a^k_{(1)}) \, d(a^k_{(2)}) =
- S(a^k_{(1)}) \, d(a^k_{(2)}) \, S(a^k_{(3)}) \, d(a^k_{(4)}) = \\
= - \chi_i(a^k_{(1)}) \chi_j(a^k_{(2)}) \, \omega^i \wedge \omega^j =
- t^k_{ij} \, \omega^i \wedge \omega^j =
- C^k_{ij} \,  \omega^i \otimes \omega^j  \; ,
\end{array}
\end{equation}
where $t^i_{mn} := < \chi_m \, \chi_n , \, a^i >$
and we use relations:
\begin{equation}
\label{23}
C^i_{jk}  =  (1 -\sigma)^{mn}_{jk} \, t^i_{mn}  \; ,
\end{equation}
which can be obtained by pairing (\ref{21}) with $a^k$.

\vspace{0.1cm}
\noindent
{\bf Remark 1.} The action of $\Delta_R$ on the first relation of
(\ref{15}) gives
$$ 
\chi_i(a_{(2)}) \omega^i \otimes S(a_{(1)})
a_{(3)}  = \chi_i(a) \omega^j \otimes r^i_j \; ,
$$ 
and as a result the
elements $r^i_j$ are expressed in terms of the generators $a^i$:
$$
\chi_i(a) \, r^i_j  =
S(a_{(1)}) \, \chi_j(a_{(2)}) \, a_{(3)}  \; \Rightarrow
r^i_j =
S(a^i_{(1)}) \, \chi_j(a^i_{(2)}) \, a^i_{(3)}  \; ,
$$
On the other hand $\chi_i(a^j \, a) = f^j_i(a)$ $\forall a \in {\cal A}$
and, thus, the elements $r^i_j$ and $f^i_j$ (which 
completely define the bicovariant
bimodule over ${\cal A}$) are fixed by elements $a^i \in {\cal A}$
and $\chi_i \in {\cal A}^*$.

\vspace{0.1cm}
\noindent
{\bf Remark 2.} For further consideration it is useful to introduce
a slightly different basis of
vector fields $\widetilde{\chi}_i$ by means of the
formula $d(a) = \omega^i \, (\widetilde{\chi}_i \triangleright
a)$ \cite{RadVl}. Comparing this formula with (\ref{15}) gives the relations
$\chi_i = f^j_i \, \widetilde{\chi}_j \; \Rightarrow \;
\widetilde{\chi}_j = S^{-1}(f^i_j) \, \chi_i$, and from
(\ref{16dd}), (\ref{16a}) we have
$$
\widetilde{\chi}_i = - S^{-1}(\chi_i) \; \Rightarrow \;
\Delta(\widetilde{\chi}_i) = \widetilde{\chi}_i \otimes 1 +
S^{-1}(f^j_i) \otimes \widetilde{\chi}_j \; .
$$
By applying $S^{-1}$ to (\ref{19}) one  obtains
\begin{equation}
\label{17c}
h \, \widetilde{\chi}_i = \widetilde{\chi}_j \, <h_{(1)} , \, S(r^j_i)
> \, h_{(2)}
\; .
\end{equation}

\subsection{Heisenberg double
$\Gamma^{\wedge} {>\!\!\!\triangleleft} \Gamma^{\wedge *}$ and Cartan calculus}

The differential algebra of Cartan calculus should be constructed
\cite{Schu}, \cite{Is}, \cite{RadVl}
as a Heisenberg double $\Gamma^{\wedge} {>\!\!\!\triangleleft} \Gamma^{\wedge *}$
of algebras $\Gamma^{\wedge}$ and $\Gamma^{\wedge *}$. 
The action of $d$ on the elements of 
the Heisenberg double $\Gamma^{\wedge} {>\!\!\!\triangleleft} \Gamma^{\wedge *}$
has to be extended as well.

First of all we recall that the Heisenberg double
$\Gamma^{\wedge} {>\!\!\!\triangleleft} \Gamma^{\wedge *}$
is an associative algebra which is a product of two algebras
$\Gamma^{\wedge}$, $\Gamma^{\wedge *}$
with nontrivial $Z_2-$ graded cross-multiplication rule
\begin{equation}
\label{25} 
\gamma \, \omega = (\gamma_{(1)} \triangleright \omega) \,
\gamma_{(2)} = (-1)^{{\rm deg} (\gamma_{(1)}) \cdot {\rm deg} (\omega_{(2)})}
\, \omega_{(1)} \, <\gamma_{(1)} , \, \omega_{(2)}> \,
\gamma_{(2)}  \; , 
\end{equation} 
where $\omega \in \Gamma^{\wedge}$ and
$\gamma \in \Gamma^{\wedge *}$. This rule 
defines the commutation relations between the elements of $\Gamma^{\wedge}$
and the elements of $\Gamma^{\wedge *}$.

Although $\Gamma^{\wedge}$ and $\Gamma^{\wedge
*}$ are Hopf algebras, their Heisenberg double  $\Gamma^{\wedge}
{>\!\!\!\triangleleft} \Gamma^{\wedge *}$ is not a Hopf algebra. But
$\Gamma^{\wedge} {>\!\!\!\triangleleft} \Gamma^{\wedge *}$ still possesses some
covariance properties. Let us define a right ${\cal A}$ - coaction and
a left ${\cal A}^*$ - coaction
on the algebra $\Gamma^{\wedge} {>\!\!\!\triangleleft} \Gamma^{\wedge *}$, which
respect the algebra structure of $\Gamma^{\wedge} {>\!\!\!\triangleleft}
\Gamma^{\wedge *}$. We denote $\{ e^{\alpha} \}$ and $\{
e_{\alpha} \}$ the dual basis elements of ${\cal A}^*$ and ${\cal
A}$ respectively. The right ${\cal A}$ -
coaction and left ${\cal A}^*$ - coaction
on $z \in \Gamma^{\wedge} {>\!\!\!\triangleleft} \Gamma^{\wedge *}$ 
are defined as follows:
\begin{equation}
\label{26} \Delta_R( z) = C \, (z \otimes 1) \, C^{-1} \; , \;\;\; 
\Delta_L( z) = C^{-1} \, (1 \otimes z) \, C \; , \;\;\;
C \equiv e^{\alpha} \otimes e_{\alpha} \; .
\end{equation}
Note that $\Delta_R(z) = \Delta(z)$ $\forall z \in {\cal A}$
and $\Delta_L(z) = \Delta(z)$ $\forall z \in {\cal A}^*$.
The axioms 
$$
(id \otimes \Delta)\Delta_R = (\Delta_R \otimes id)\Delta_R
\;\; , \;\;\;
(id \otimes \Delta_L)\Delta_L = (\Delta \otimes id)\Delta_L \; ,
$$ 
$$
(id \otimes \Delta_R)\Delta_L(z) = C^{-1}_{13} \,
(\Delta_L \otimes id)\Delta_R(z) \, C_{13} \; , 
$$ 
can be verified directly by using the pentagon identity \cite{BSk} for $C$
$$
C_{12} \, C_{13} \, C_{23} = C_{23} \, C_{12} \; .
$$
It is clear that the coactions (\ref{26}) are
covariant transformations (homomorphisms) of the algebra 
$\Gamma^{\wedge} {>\!\!\!\triangleleft}
\Gamma^{\wedge *}$. The inverse of the canonical element $C$ is
$$
C^{-1} = S(e^{\alpha}) \otimes e_{\alpha} = e^{\alpha} \otimes S(e_{\alpha}) \; ,
$$ 
and $\Delta_R$ (\ref{26}) is rewritten in the form
\begin{equation}
\label{27}
\Delta_R( z) =
(e^{\gamma}_{(1)} \, z \, S(e^{\gamma}_{(2)})) \otimes e_{\gamma} \; .
\end{equation}
In particular, for $y \in {\cal A}^*$ we have
\begin{equation}
\label{28}
\Delta_R (y) = y_{(1)} \otimes y'_{(2)} =
(e^{\gamma}_{(1)} \, y \, S(e^{\gamma}_{(2)})) \otimes e_{\gamma} \; .
\end{equation}
By pairing the second factor of (\ref{28}) with an arbitrary $\;\; x \in
{\cal A}^* \;\;$ we deduce the relation $y_{(1)} \, < x, \, y'_{(2)} > = x_{(1)} \,
y \, S(x_{(2)})$ which is equivalent to the commutation relations
for $x,y \in {\cal A}^*$
\begin{equation}
\label{29}
x \, y = y_{(1)} \, < x_{(1)}, \, y'_{(2)} > \, x_{(2)} \; ,
\end{equation}
where $\Delta(x) = x_{(1)} \otimes x_{(2)}$. The inverse statement
(that one can obtain (\ref{28}) from (\ref{29})) is also correct.
The analogs of eqs. (\ref{27}), (\ref{28}) and (\ref{29}) 
for coaction $\Delta_L$ can be deduced in the same way.
Comparing the commutation relations (\ref{29}) 
with eqs. (\ref{17c}) one can find as an example the right coaction
$\Delta_R$ on $\widetilde{\chi}_i$ (\ref{28}):
\begin{equation}
\label{31a}
\Delta_R(\widetilde{\chi}_i) = \widetilde{\chi}_j \otimes S(r^j_i) \; .
\end{equation}

At the end of this Section we present
some cross-commutation relations (see (\ref{25}))
which will be needed below
\begin{equation}
\label{com4}
[\chi_i , \, \omega^j ]_- =  C^j_{lk} \omega^l \, f^k_i  \;
, \;\;\;
[\gamma_i , \, \omega^j]_+ = f^j_i  , \;\;\;
[\gamma_i , \, a] = 0 \;\; \forall a \in {\cal A} \; .
\end{equation}

\section{BRST Operator for Quantum Lie Algebras}

In this Section we find a bi-invariant element $Q \in
\Gamma^{\wedge} {>\!\!\!\triangleleft} \Gamma^{\wedge *}$
(a BRST operator) which generates the differential $d$: 
\begin{equation}
\label{Qdiff}
d\omega = [ Q, \, \omega ]_{\pm} \; , \;\;\; \forall \omega \in \Gamma^{\wedge}
\end{equation}
The operator $Q$ has to be of the grading $1$ and 
obeys $Q \wedge Q =0$.

We change the basis of differential forms $\{
\omega^i \}$  and consider new basics elements
$\Omega^i = \omega^j \, S^{-1}(f^i_j)$ which are convenient
in view of the relation
$$
[ \Omega^i , \, a ]  = 0 \;\;\;\; \forall a \in
{\cal A} \; .
$$
The following equations for $\Omega^i$ are also valid (see (\ref{com4})):
\begin{equation}
\label{app2.1} 
\chi_{|2>} \, \Omega^{<2|} =  \Omega^{<1|} \, \left(
\sigma_{12} \, \chi_{|1>} +  C^{<2|}_{|12>} \right) \, , 
\;\; \gamma_{|2>} \, \Omega^{<2|} = - \Omega^{<1|} \,
\sigma^{-1}_{12} \, \gamma_{|1>} + I_2  \; ,
\end{equation}
and the definition of the
wedge product for $r$ variables $\Omega^i$ holds:
$$
\Omega^{<r|} \dots \Omega^{<1|} = \Omega^{<r|} \otimes \dots
\otimes \Omega^{<1|} \, A_{1 \rightarrow r} \; .
 $$

\vspace{0.1cm}
\noindent
Now we formulate the main result of the paper
(see also \cite{BIO}): 
 
 \vspace{0.1cm}
 \noindent
 {\bf Proposition}
 {\it The BRST operator $Q$ for the quantum algebra (\ref{int1}),
 which generate the differential (\ref{Qdiff}), has the following form
  \begin{equation}
 \label{brst} Q = \Omega^i \, \chi_i + \sum^{h-1}_{r=1} \, Q_{(r)} \; ,
  \end{equation}
 where $h$ is the height of the operator $\sigma_{12}$ (\ref{anti}).
 
 Here the operators $Q_{(r)}$ are given by
  \begin{equation}
 \label{qrr} Q_{(r)} =\Omega^{<r+1|} \, \Omega^{<r|} \dots
 \Omega^{<1|} \, X^{<\tilde{1} \dots \tilde{r}|}_{|1 \dots r+1>} \,
 \gamma_{|\tilde{1}>} \dots \gamma_{|\tilde{r}>}
  \end{equation}
 (the wedge product is implied); $X^{<1 \dots
 r |}_{|1 \dots r+1>}$ are tensors which satisfy the following recurrent
 relations
  \begin{equation}
  \label{app3.3}
 A_{1 \rightarrow r+1} \,
  X^{< 1\dots r|}_{|1 \dots r+1>} \, A_{1 \rightarrow r } = 
 A_{1 \rightarrow r+1} \, \left((-1)^{r} \, 
 \sigma_{r+ 1 \leftarrow 1} - {\bf 1} \right) \, 
 X^{< 2\dots r |}_{|2 \dots r+1>} \,
 A_{2 \rightarrow
  r}  \; 
  \end{equation}
 with the initial condition $A_{12} \, X^{<0|}_{|12>} = - C^{<0|}_{|12>}$.
  }
  
\vspace{.1cm}
\noindent
{\bf Proof}. We have to verify the conditions
$[Q, \, a]=d \, a$ $(\forall a \in {\cal A})$ and $[Q, \,
\omega^i]_+ =d \, \omega^i$ $(\forall \omega^i \in \Gamma^{(1)})$
and the identity $Q^2 = 0$.
The proof of the condition $[Q, \, a]=d \, a$ $(\forall a \in {\cal
A})$ is straightforward since $\Omega^i$ and $\gamma_i$ commute with all $a \in
{\cal A}$. So we need only to prove $[\Omega^i \, \chi_i , \, a]=d \,
a$ which follows immediately from relations (\ref{15}).
The other conditions are valid only if the recurrent relations
(\ref{app3.3}) are fulfilled. We will give the complete proof of this
statement elsewhere.

\vspace{0.1cm}
 \noindent
 {\bf Remark 1.}  The first term  $\Omega^i \,
\chi_i = \omega^i \, \widetilde{\chi}_i$ in (\ref{brst})
is a bi-invariant element (see (\ref{2}), (\ref{31a})). 
One can prove that all other terms $Q_{(r)}$ are also
bi-invariants.

 \vspace{0.1cm}
\noindent
 {\bf Remark 2.}
 For general $\sigma^{ij}_{kl}$ and $C^i_{jk}$
  it is rather difficult to solve equations (\ref{app3.3}) explicitly. 
  However for the case $\sigma^2 = {\bf 1}$ the main equations 
  (\ref{app3.3}) become simpler and the general solution for $Q$
  can be found. Indeed the
  relation (\ref{app3.3}) for $r=2$ gives
 \begin{equation}
 \label{r=2}
  A_{1 \rightarrow 3} \,
  X^{<12|}_{|123>} \, \left({\bf 1} - \sigma_{12}  \right) = A_{1
  \rightarrow 3} \, \left( \sigma_{23} \, \sigma_{12} -{\bf 1} \right) \,
  X^{<2|}_{|2 3>}  \; .
\end{equation}
 For $\sigma^2 =1$ 
 we have $A_{1 \rightarrow 3} \, 
 \left( \sigma_{23} \, \sigma_{12} -{\bf 1} \right) = 0$ and 
 therefore $Q_{(r)} = 0$ for $r \geq 2$. Thus the BRST operator
 (\ref{brst}) has the familiar form
 $$
 Q=\Omega^{<1|} \, \chi_{|1>} - \Omega^{<2|} \otimes \Omega^{<1|}  \,
 C^{<1|}_{|1 2>} \, \gamma_{|1>} \; .
 $$ 
 
 In general, for $\sigma^2 \neq 1$, the sum in (\ref{brst})
 will be limited only by the height $h$
 of the operator $\sigma$.
 
 Below we present an explicit form for $Q$ 
 for the standard quantum deformation
 ${\cal A}^* = U_q(gl(N))$ of the universal enveloping algebra
 of the Lie algebra $gl(N)$ ($\sigma^2 \neq 1$ in this case).

\vspace{0.1cm}
\noindent
{\bf Remark 3.} Here we give expressions for first two coefficients
$X^{<12|}_{|123>}$ and $X^{<123|}_{|1234>}$. Substitution of
the initial condition $A_{12} \, X^{<0|}_{|12>} = - C^{<0|}_{|12>}$
into (\ref{r=2}) gives
\begin{equation}
\label{x1}
A_{123} \, X^{<12|}_{|123>} \,  A_{12} = 
 \left[ C_{2}  + \sigma_{1}  \,  \sigma_{2} \, C_{1} 
 \, \delta_{3} \right] \,  A_{12}\; 
\end{equation}
where we have used the concise notation
\begin{equation}
\label{nota}
\sigma_n = \sigma_{n \, n+1} \; , \;\;\; C_{n} = C^{<n|}_{|n \, n+1>} \; , 
\;\;\; \delta_{n} = \delta^{<n-1|}_{|n>} \; .
\end{equation}
The analogous formula for the next coefficient is
\begin{equation}
\label{x2}\begin{array}{l}
 A_{1234} \, X^{<123|}_{|1234>} \,  A_{123}\\ = 
 \left[ (C_{3}  + \sigma_{2}  \,  \sigma_{3} \, C_{2} \, \delta_4) 
- \sigma_{1}  \,  \sigma_{2}   \,  \sigma_{3} \, 
(C_{2}  + \sigma_{1}  \,  \sigma_{2} \, C_{1} \, \delta_3) \delta_4
 \right] \,  A_{123} \; .
\end{array}\end{equation} 
To obtain (\ref{x1}) and (\ref{x2}) 
it is convenient to rewrite eqs. (\ref{int2}) - (\ref{int1a})
in the form (see notation (\ref{nota})):
$$
 \begin{array}{c}
C_{1} \, \delta_3 \, C_{1} = \sigma_{2} \, C_{1} \,
\delta_3 \, C_{1} + C_{2} \, C_{1} \; , \;\;\;
C_{1} \, \delta_3 \, \sigma_{1} = 
\sigma_{2} \, \sigma_{1} \, C_{2} \; , \\ \\
(\sigma_{2} \, C_{1} \, \delta_3 +  C_{2} ) \, \sigma_{1} = 
\sigma_{1} \, (\sigma_{2} \, C_{1} \, \delta_3 + C_{2} ) \; , \;\; \;
 \sigma_{1} \, \sigma_{2} \, \sigma_{1} =
\sigma_{2} \, \sigma_{1} \, \sigma_{2} \; .
\end{array}
$$ 

\section{BRST and anti-BRST operators 
for quantum linear algebra $U_q(gl(N))$.}

\subsection{BRST operator for $U_q(gl(N))$}

The quantum algebra $U_q(gl(N))$
is defined (as a Hopf algebra)
by the relations \cite{FRT}
\begin{equation}
\label{92}
\hat{R} \, L^{\pm}_2 L^{\pm}_1  =  L^{\pm}_2 L^{\pm}_1 \, \hat{R} \,, \ \ \  
\hat{R} \, L^+_2 L^-_1  = L_2^- L^+_1 \, \hat{R} \,,  
\end{equation}
\begin{equation}
\label{93}
\Delta(L^{\pm})  =  L^{\pm}\otimes L^{\pm} , \ \ \
\varepsilon(L^{\pm}) = {\bf 1} , \ \ \
S(L^{\pm}) = (L^{\pm})^{-1} \; ,
\end{equation}
where elements of the $N \times N$
matrices $(L^\pm)^i_j$ are generators of $U_q(gl(N))$;
the matrices $L^+$ and $L^-$ are 
respectively upper and lower triangular, 
their diagonal elements are related by $(L^+)^i_i \, (L^-)^i_i = 1$
for all $i$.
The matrix $\hat{R}$ is defined as
$\hat{R} := \hat{R}_{12}= P_{12} \, R_{12}$ ($P_{12}$ is the permutation matrix) and
the matrix $R_{12}$ is the standard
Drinfeld-Jimbo $R$-matrix for $GL_q(N)$, 
\begin{equation}
\label{Rgl}
R_{12}=R^{i_{1},i_{2}}_{j_{1},j_{2}}=
\delta^{i_{1}}_{j_{1}} \delta^{i_{2}}_{j_{2}}(1+(q-1)\delta^{i_{1}i_{2}}) +
(q-q^{-1})\delta^{i_{1}}_{j_{2}} \delta^{i_{2}}_{j_{1}}
\Theta_{i_{1}i_{2}} \; , 
\end{equation}
where
$$
\Theta_{ij} = 
\left\{ 
\begin{array}{c}
1 \;\;\; {\rm if} \;\; i>j   \; ,  \\
0 \;\;\; {\rm if} \;\; i \leq j  \; .
\end{array}
\right.
$$
This $R$-matrix satisfies the Hecke condition 
$\hat{R}^2 = \lambda \, \hat{R} + {\bf 1}$, where $\lambda = (q-q^{-1})$
and $q$ is a parameter of deformation.

The generators of the algebra ${\cal A}^*$ are defined by the formula  
\cite{Jur}, \cite{Is}, \cite{RadVl}
\begin{equation}
\label{102} 
\chi _{k}^{l}
=\frac{1}{\lambda}\,[\,(D^{-1})_{k}^{l}-(D^{-1})_{i}^{j}f_{kj}^{li}\,] \; .
\end{equation}
Here 
$f_{kj}^{li}={L^-}_{k}^{i} S({L^+}_{j}^{l})$ and the
numerical matrix $D$ can be found by means of relations
$$
Tr_2 \hat{R}_{12} \Psi_{23} = P_{13} = Tr_2  \Psi_{12} \, \hat{R}_{23}  \; , \;\;\;
D_1 := Tr_2 \Psi_{12} \; \Rightarrow \; Tr_1(D^{-1}_1 \hat{R}^{-1}) = {\bf 1}_2 \; ,
$$
where $Tr_1$ and $Tr_2$ denote the traces over first and second spaces.

For the $GL_q(N)$ $R$-matrix (\ref{Rgl}) the
explicit expression for the $D$-matrix is:
\begin{equation}
\label{D}
(D^{-1})^i_j = q^{2(N-i) +1} \, \delta^i_j  , \;\;
Tr (D^{-1}  ) = \frac{q^{2N} - 1}{q-q^{-1}}  .
\end{equation}

It is convenient to write down
the commutation relations for the differential algebra
in terms of generators 
$$
L^i_j = (L^+)^i_k S((L^-)^k_j) = \delta^i_j - \lambda \, 
S^{-1}(\chi^i_k) \, D^k_j \; , 
$$
$$
J^i_n = - S^{-1}(f^{ik}_{jl}) \, \gamma^l_k D^j_n \; , \;\;\;
\omega^i_j = \Omega^k_m \, f^{mi}_{kj} \; .
$$
The indices now are pairs of indices;
the roles of the elements $\chi_i$, $\gamma_j$ and $\Omega^k$
are played by the generators $\chi^i_j$, $\gamma^i_j$ and $\Omega^i_j$
respectively.

The commutation relations are \cite{SWZ}, \cite{Is}, \cite{RadVl}:
\begin{equation}
\label{106} 
\omega_2 \hat{R}^{-1} \omega_2 \hat{R} = - \hat{R}^{-1} \omega_2 \hat{R}^{-1} \omega_2
\; , \;\;\;
\omega_2 \, \hat{R} \, L_2 \, \hat{R} = \hat{R} \, L_2 \, \hat{R} \, \omega_2 \; , 
\end{equation}
\begin{equation}
\label{109}
\omega_2 \, \hat{R} \, J_2 \hat{R} + \hat{R} \, J_2 \, \hat{R} \, \omega_2 = -\hat{R} \, ,
\;\;\;
L_2  \hat{R} \, L_2  \hat{R} = \hat{R} \, L_2 \hat{R} \, L_2 \; , 
\end{equation}
\begin{equation}
\label{113}
J_2 \hat{R} \, L_2 \hat{R} = \hat{R} \, L_2 \hat{R} \, J_2 \; , \;\;\;
J_2 \hat{R} \, J_2 \hat{R} = - \hat{R}^{-1} J_2 \hat{R} \, J_2 \, .
\end{equation}

To write down the
whole differential algebra over $GL_q(N)$ we need to add the
generators $T^i_j$ of the quantum group $Fun(GL_q(N))$ 
with commutation relations
\begin{eqnarray}
\label{115}
 \hat{R} \, T_1 \, T_2 = T_1 \, T_2 \, \hat{R} \; , \;\;
\omega_1 \, T_2  = T_2 \, \hat{R}^{-1} \, \omega_2 \, \hat{R}^{-1}  ,
\end{eqnarray}
\begin{eqnarray}
\label{113a}
J_1 \, T_2  = T_2 \, \hat{R} \, J_2 \, \hat{R} \; , \;\;
L_1 \, T_2  = T_2 \, \hat{R} \, L_2 \, \hat{R}  \; .
\end{eqnarray}

The BRST operator $Q$ for the differential algebra 
(\ref{106}) -- (\ref{113a}) can be constructed with the help
of the formula (\ref{brst}) and has the following form \cite{BIO}:
\begin{equation}
\label{Qgl2}
Q = Tr_q \left( \omega \,(L- {\bf 1})/\lambda 
-    \omega \, L \, (\omega J) + \lambda \, \omega \, L \, (\omega J)^2 -
\lambda^2 \, \omega \, L \, (\omega J)^3 + \dots \right)
\end{equation}
\begin{equation}
\label{Qgl}
= Tr_q \left( \omega \, \frac{(L- {\bf 1})}{\lambda} 
- \omega \, L \, (\omega J) \, ({\bf 1} + \lambda \omega J)^{-1} \right)
= - \frac{1}{\lambda} \, Tr_q ( \omega ) + \frac{1}{\lambda} \,
Tr_q \left( \Theta  \right)
\; ,
\end{equation}
where $\Theta = \omega \, L\, ({\bf 1} + \lambda \, \omega J)^{-1}$
and $Tr_q(Y) := Tr(D^{-1} \, Y)$ is a quantum trace.
The sum in eq. (\ref{Qgl2})
is finite due to the fact that monomials in $\omega$'s
of the order ${N^2+1}$ are equal to zero (since the exterior algebra 
(\ref{106}) of forms on the quantum group $GL_q(N)$
is a flat deformation of the classical algebra 
\cite{Man}, \cite{IsPya}).

One can check directly that the BRST operator $Q$
given by (\ref{Qgl}) satisfies:
\begin{equation}
\label{diffq}
 Q^2 = 0 \; ,
\;\;\; [ Q , \, L ] = 0 \; ,
\end{equation}
\begin{equation}
\label{diffT}
[Q , \, T] = T \, \omega \equiv d \, T , \;\;
[Q , \, \omega]_+ = - \omega^2 \equiv d \, \omega ,
\end{equation}
\begin{equation}
\label{cartan}
[ Q , \, J ]_+ = \frac{1}{\lambda} \, ({\bf 1}- L)  \; .
\end{equation}
The (anti)commutator with the BRST operator $Q$
(relations (\ref{diffT})) defines the exterior differential operator
over $GL_q(N)$; 
it provides the structure of the de Rham complex $\Omega(GL_q(N))$ on the
subalgebra with generators $T^i_j$ and $\omega^i_j$
(the de Rham complex $\Omega(GL_q(N))$ has been firstly considered
by Yu.I.Manin, G.Maltsiniotis and B.Tsygan \cite{Man}).

The last relation (\ref{cartan}) is an analog of the Cartan identity.
To obtain relations (\ref{diffq}) -- (\ref{cartan})
one has to use the invariance property
of the quantum trace:
$$
Tr_q(X){\bf 1}_2 = Tr_{q1}(\hat{R}^{\pm 1} X_2 \hat{R}^{\mp 1})
$$
and relations
\begin{equation}
\label{th}
\hat{R} \, \Theta_2 \, \hat{R}^{-1} \, \omega_2
= - \omega_2 \, \hat{R}^{-1} \, \Theta_2 \, \hat{R}
\; ,
\end{equation}
\begin{equation}
\label{th1}
\hat{R} \, \Theta_2 \, \hat{R}^{-1} \, \Theta_2 =
- \Theta_2 \, \hat{R}^{-1} \, \Theta_2 \, \hat{R}^{-1}  \; ,
\end{equation}
\begin{equation}
\label{TL}
\hat{R}^{-1} \, \Theta_2 \, \hat{R} \, L_2 =
L_2 \, \hat{R} \, \Theta_2 \, \hat{R}^{-1}  \; , \;\;\;
 \Theta_1 \, T_2 =
T_2 \, \hat{R}^{-1} \, \Theta_2 \, \hat{R}  \; ,
\end{equation}
\begin{equation}
J_2 \, \hat{R} \, \Theta_2 \, \hat{R}^{-1} +
\hat{R}^{-1} \, \Theta_2 \, \hat{R} \, J_2  
= - L_2 \, ({\bf 1} + \lambda \, \omega J)^{-1}_2 \, \hat{R}^{-1} \,
({\bf 1} + \lambda \, \omega J)_2 \; , \label{JT}
\end{equation}
which follow from  (\ref{106})-(\ref{113a}).

\vspace{0.1cm}
\noindent
{\bf Remark 1.} The operator $Q$ given by (\ref{Qgl}) has the
correct classical limit for $q \rightarrow 1$ ($\lambda \rightarrow 0$,
$L \rightarrow {\bf 1} + \lambda \tilde{\chi}$,
$\omega \rightarrow \tilde{\omega}$,
$J \rightarrow - \tilde{\gamma}$):
$$
Q \rightarrow Q_{cl}=
Tr(\tilde{\omega} \tilde{\chi} + \tilde{\omega}^2 \, \tilde{\gamma} ) =
Tr(\tilde{\omega} \, X - \tilde{\omega} \, \tilde{\gamma} \, \tilde{\omega})
\; ,
$$
where 
\begin{equation}
X := \tilde{\chi} + \tilde{\omega}
\, \tilde{\gamma} + \tilde{\gamma} \, \tilde{\omega}
\label{defofx}\end{equation}
 and
the classical algebra is
\begin{equation}
[\tilde{\omega}_2, \, \tilde{\gamma_1} ]_+ = P_{12} \; , \;\;\;
 [\tilde{\omega}_2, \, \tilde{\omega}_1 ]_+ = 0
= [\tilde{\gamma}_2, \, \tilde{\gamma}_1 ]_+ \; ,
\end{equation}
\begin{equation}
\label{xx}
[ X_2 , \, X_1 ] = P_{12}(X_2 - X_1) \; , 
\end{equation}
\begin{equation}
\label{xxx}
 [X_2, \, \tilde{\omega}_1 ] =  0  = [X_2, \, \tilde{\gamma}_1 ] \; .
\end{equation}

\vspace{0.1cm}
\noindent
{\bf Remark 2.} The differential complex $\Omega(GL_q(N))$ 
gives rise to the de Rham cohomology groups $H^p(GL_q(N))$ 
of the quantum groups $GL_q(N)$. 
The $q$-analogs of the basic generators for the de Rham cohomology ring
$H^*(GL_q(N))$ can be chosen as
\begin{equation}
\Omega^{(n)} := Tr_q( \omega^n ) \;\;\; (n=1,3,5, \dots ,2N-1) \; .
\end{equation}
These generators satisfy \cite{Is}, \cite{FP}
$$
\, d \, \Omega^{(n)} = -Tr_q( \omega^{n+1} ) = 0  , \;\;
\, [\Omega^{(n)}, \, \Omega^{(k)}]_+ = 0 \; . 
$$

\subsection{Anti-BRST operator and quantum
Laplacian}

In the same way
as we deduce the explicit formula (\ref{Qgl2})-(\ref{Qgl})
for the BRST operator $Q$, 
one can construct the anti-BRST operator $Q^*$
for the algebra $U_q(gl(N))$:
\begin{equation}
\label{Q*}
Q^* =  Tr_q \left( J \, (L^{-1} - {\bf 1})/ \lambda
+ J \, L^{-1} \, J \, \omega \right)
=  \frac{1}{\lambda} \, \left( \, Tr_q ( \Theta^*) - Tr_q ( J )  \, \right)
\; ,
\end{equation}
where $\Theta^* = J \, L^{-1} \, ({\bf 1} + \lambda J \omega)$.
The operator $Q^*$ satisfies
\begin{equation}
\label{Q*L}
(Q^*)^2 = 0 \; , \;\;\; [ Q^* , \, L ] = 0 \; ,  \;\;\;
\end{equation}
\begin{equation}
\label{Q*T}
[ Q^* , \, T ] = q^{2N} \, T \, J , \;\;
[Q^* , \, J ]_+ = - q^{2N} \, J^2  \; ,
\end{equation}
\begin{equation}
\label{Q*om}
[ Q^* , \, \omega ] = q^{2N} \left( W \, \frac{1 - L^{-1}}{\lambda} \,
\overline{W} - \lambda \, \omega \, J^2 \, \omega \right) .
\end{equation}
Here the factor $q^{2N}$ appeared because of
the following identities (see (\ref{D}))
$$
Tr_{q1}(R) = Tr_{q1}(R^{-1} + \lambda{\bf 1}) =
q^{2N} \, {\bf 1}_2 \; .
$$
For convenience we introduce new operators \cite{Is}
\begin{equation}
\label{wbarw}
W = ({\bf 1} + \lambda \, \omega \, J) \; , \;\;\;
\overline{W} = ({\bf 1} + \lambda \, J \, \omega ) \; .
\end{equation}
To prove eqs. (\ref{Q*L}), (\ref{Q*T}) and (\ref{Q*om})
we used relations
\begin{equation}
\label{T*com}
\begin{array}{l}
\Theta^*_1 \, T_2  =  T_2 \, \hat{R} \, \Theta^*_2 \, \hat{R}^{-1} \; ,
\;\;\;
\hat{R}^{-1} \, \Theta^*_2 \, \hat{R} \, J_2 = 
- J_2 \, \hat{R} \, \Theta^*_2 \, \hat{R}^{-1}
\; , \\[1em]
\hat{R} \, \Theta^*_2 \, \hat{R} \, \Theta^*_2 = 
- \Theta^*_2 \, \hat{R} \, \Theta^*_2 \, \hat{R}^{-1}  \; ,
\;\;\;
\hat{R}^{-1} \, \Theta^*_2 \, \hat{R} \, L_2  = 
L_2 \, \hat{R} \, \Theta^*_2 \, \hat{R}^{-1} \; , \\[1em]
\hat{R} \, \Theta^*_2 \, \hat{R}^{-1} \, \omega_2  + 
\omega_2 \, \hat{R}^{-1} \, \Theta^*_2 \, \hat{R}  =
- (W \, L^{-1} \, \overline{W})_2 \, \hat{R} \; .
\end{array}
\end{equation}

Using eqs. (\ref{Q*T}) and (\ref{Q*om}) 
one can define the dual differential $d^*$
as an (anti)commutator with the anti-BRST operator $Q^*$.

Now we define the current matrix $U^i_j$:
$$
U :=  W \, L^{-1} \, \overline{W} = 
- q^{-2N} \, [Tr_q(\Theta^*), \, \omega]_+ \; ,
$$ 
which is a quantum analog of the
matrix $X$ (\ref{defofx}) and satisfies the reflection equation
$$
\hat{R}^{-1} \, U_2 \, \hat{R}^{-1} \, U_2 = 
U_2 \, \hat{R}^{-1} \, U_2 \, \hat{R}^{-1} \; .
$$
For this current matrix we deduce the following commutation relations
\begin{equation}
\label{uo}
\hat{R} \, U_2 \, \hat{R}^{-1} \, \omega_2
 =  \omega_2 \, \hat{R}^{-1} \, U_2 \, \hat{R}  \; ,  \;\;\;
\hat{R}^{-1} \, U_2 \, \hat{R} \, J_2
 =   J_2 \, \hat{R} \, U_2 \, \hat{R}^{-1} \; . 
\end{equation}
These relations are correct quantum analogs of 
the commutativity conditions (\ref{xxx})
($U$ is represented as $U = {\bf 1} + \lambda {\bf X}$ where
${\bf X} \rightarrow X$ for $q \rightarrow 1$),
i.e. 
$$
 {\bf X}_{|1>} \, \omega_{|2>} = 
\sigma_{12} \, \omega_{|2>} \, {\bf X}_{|1>}  \; , \;\;\;
 {\bf X}_{|1>} \, J_{|2>} = 
\overline{\sigma}_{12} \, J_{|2>} \, {\bf X}_{|1>} \; ,
$$ 
where $\sigma_{12}$, $\overline{\sigma}_{12}$ are braid matrices defined by
(\ref{uo}) converted to the permutation
matrices for $q = 1$. Moreover,
the transformation of the left-invariant current $U$ into the right invariant
current $\Xi$:
$$
\Xi = T \, W \, L^{-1} \, \overline{W} \, T^{-1}
= T \, U \, T^{-1} \; ,
$$
leads to the exact commutativity conditions:
$[\Xi_2 , \, \omega_1] = 0$, $[\Xi_2, \, J_1] = 0$.

By definition the quantum Laplace operator is
(here it is enough to use relations (\ref{diffq}) - (\ref{cartan}))
\begin{equation}
\label{Lapl}
\Delta := Q \, Q^* + Q^* \, Q 
= \lambda^{-2} \, Tr_q \left( L  - 2  +
L^{-1} \, \overline{W}  - \lambda \,
J \, L^{-1} \, \overline{W} \, \omega \right) , 
\end{equation}
The Laplacian $\Delta$ is a BRST and anti-BRST invariant operator:
$$
[Q, \, \Delta ] = 0 = [Q^*, \, \Delta] \; ,
$$
and it generalizes the Casimir operator for the universal enveloping
algebra $U_q(gl(N))$.

Taking into account the identity
$$
- \lambda \,
Tr_q \left( J \, L^{-1} \, \overline{W} \, \omega \right) 
= Tr_{q} \left( (q^{2N} - 1) \, L^{-1} \, \overline{W} + \lambda \,
q^{2N} \, \omega \, J \, L^{-1} \, \overline{W} \right) \; ,
$$
one obtains a remarkable expression for the quantum Laplacian via the
current $U$:
$$
\Delta = \frac{1}{\lambda^2} \, Tr_q \left( L  +
q^{2N} \, U - 2 \right) \; .
$$

Now the formulation of the Hodge decomposition theorem
is in order.
Consider the space of polynomials in the variables $\omega^i_j$ and
$T^i_j$ with complex coefficients $\psi^{i_1 \dots i_k}_{j_1 \dots j_k}$,
$$	
|\Psi>  :=  \Psi[T,\omega] 
=  {\displaystyle \sum_{k=0}^K \sum_{r=0}^{k}}Tr_{1 \dots k} 
\left(T_1 \dots T_r \, \omega_{r+1} \dots \omega_k \psi_{1 \dots k} \right) 
\; 
$$
(for some $K$). Here the case $r=0$ corresponds to the arbitrary polinomial 
in $\omega$'s which independent of $T$'s.
The vector fields $L^i_j -\delta^i_j$ and inner derivatives $J^i_j$
act on the zero-order monomial $|0> :=1$ from the left as anihilation operators:
$$
(L^i_j -\delta^i_j) \, |0> = 0 \; , \;\;\; J^i_j \, |0> = 0 \; .
$$
This defines the left action of the BRST, anti-BRST and Laplace operators on the
polynomials $|\Psi>$. Now the decomposition theorem can be formulated: \\

\noindent
{\bf Theorem.} {\it Any polynomial $\Psi[T,\omega]$ can be decomposed
into a sum of BRST-exact, co-exact and harmonic polynomials:
$$
|\Psi> = |\Omega> + Q \cdot |\chi> + Q^* \cdot |\Phi> \; ,
$$
where $\Delta \, |\Omega> = 0$.
}

The proof of this theorem is straightforward and analogous to the
proof of the decomposition theorem
in the case of the classical Lie algebras (see e.g. \cite{vanH}, \cite{CAMP}).

\vspace{0.1cm}
\noindent
{\bf Acknowledgements.} This work was partially supported
by RFBR grants 98-01-22033 and the CNRS grant PICS-608. The work of
AI was also supported by RFBR grant 00-01-00299.

\end{document}